\documentclass[11pt]{article}
\usepackage{amsfonts,latexsym,amssymb,epsfig}
\usepackage{amsmath,amsthm,amstext}
\usepackage{fullpage}
\newcommand{\shlomo}[1]{}
\newcommand{\avner}[1]{}

\begin{document}
\bibliographystyle{plain}
\newtheorem{theorem}{Theorem}
\newtheorem{proposition}[theorem]{Proposition}
\newtheorem{corollary}[theorem]{Corollary}
\newtheorem{note}[theorem]{Note}
\newtheorem{lemma}[theorem]{Lemma}
\newtheorem{definition}[theorem]{Definition}
\newtheorem{conjecture}[theorem]{Conjecture}
\newcounter{fignum}
\newcommand{\figlabel}[1]
           {\\Figure \refstepcounter{fignum}\arabic{fignum}\label{#1}}
\newcommand{\ignore}[1]{}
\font\boldsets=msbm10
\def\Z{{\hbox{\boldsets  Z}}}
\def\R{{\hbox{\boldsets  R}}}
\def\half{\frac 1 2}
\def\F2{{\{0,1\}}}
\def\A{{\cal A}}
\def\eps{{\epsilon}}
\def \mixx {\mbox {mix}}
\newcommand{\rank} {\mbox {rank}}
\newcommand{\schreier} {\mbox {sc}}
\newcommand{\Hom} {\mbox{Hom}}
\newcommand{\G} {{\cal F}}
\newcommand{\gr}[9]{
  \fbox{$\begin{array}{ccc}#1&#2&#3\\#4&#5&#6\\#7&#8&#9\end{array}$}
}
\title{A counterexample to a conjecture of Bj\"{o}rner and Lov\'asz on the 
$\chi$-coloring complex}
\author{
Shlomo Hoory
\\Department of Computer Science
\\University of British Columbia
\\shlomoh@cs.ubc.ca
\and
Nathan Linial
\\Department of Computer Science
\\Hebrew University
\\nati@cs.huji.ac.il
}

\maketitle

\def\Hom{\mbox{Hom}}

\begin{abstract}
Associated with every
graph $G$ of chromatic number $\chi$ is another graph $G'$.
The vertex set of $G'$ consists of all $\chi$-colorings of $G$,
and two $\chi$-colorings are adjacent when they differ on exactly one vertex.
According to a conjecture of Bj\"{o}rner and Lov\'asz, this graph $G'$ must
be disconnected. In this note we give a counterexample to this conjecture.
\end{abstract}


One of the most disturbing problems in graph theory is that we only have
few methods to prove lower bounds on the chromatic number of 
graphs. A famous exception is Lov\'asz's~\cite{Lo78} proof
of the Kneser conjecture. This paper has indeed introduced a new method
into this area and
is one of the first applications of topological methods to combinatorics.
It shows how to use
the Borsuk-Ulam Theorem to derive a (tight) lower bound for the
chromatic number of the Kneser graph.
Since then, the idea of finding topological obstructions to graph colorings 
has been extensively studied~\cite{BaKo04, Ma03}. 
In particular, Bj\"{o}rner and Lov\'asz made a conjecture generalizing the 
concept of a topological obstruction to graph coloring 
(see \cite{BaKo03} conjecture 1.6).
In this note we provide a counterexample to this general conjecture.

To state the general conjecture, we need some definitions.
For two graphs $G$, $H$, an $H$-coloring of $G$ is a homomorphism from 
$G$ to $H$.
Namely, a mapping $\phi:V(G) \rightarrow V(H)$ such that for all edges 
$(x,y) \in E(G)$ we have $(\phi(x),\phi(y))\in E(H)$. 
The coloring complex $\Hom(G,H)$ is a CW-complex whose $0$-cells are the
$H$-colorings of $G$. The cells of $\Hom(G,H)$ are the maps 
$\eta$ from the vertices of $G$ to non-empty vertex subsets of $H$ such that
$\eta(x) \times \eta(y) \subseteq E(H)$ for all $(x,y) \in E(G)$.
The closure of a cell $\eta$ consists of all cells $\tilde{\eta}$ such 
that $\tilde{\eta}(v) \subseteq \eta(v)$ for all $v \in V(G)$.
We say that a complex $C$ is $k$-connected if every map from $S^k$ to $C$ 
can be extended to a map from $B^{k+1}$ to $C$. Equivalently, if all the 
homotopy groups up to dimension $k$ are trivial. 
Specifically, $(-1)$-connected means non-empty, and $0$-connected is connected.
We denote the chromatic number of a graph $G$ by $\chi(G)$.
In an attempt to capture some topological obstructions
to low chromatic number,
Bj\"{o}rner and Lov\'asz have made the following conjecture:

\begin{conjecture} \label{conj} (Bj\"{o}rner and Lov\'asz)
Let $G,H$ be two graphs such that the coloring complex $\Hom(G,H)$ is
$k$-connected. Then $\chi(H) \geq \chi(G)+k+1$.
\end{conjecture}

There are three special cases of this conjecture that are known to be true.
One is essentially
Lov\'asz's original argument implying the Kneser conjecture. The
other two are recent results of Babson and Kozlov:
\begin{theorem} 
Conjecture~\ref{conj} holds when (i) $G=K_2$ (Lov\'asz~\cite{Lo78}),
(ii) $G=K_m$ (Babson and Kozlov~\cite{BaKo03}),
(iii) $G= C_{2r+1}$ (Babson and Kozlov~\cite{BaKo04}).
\end{theorem}

We refute the conjecture by exhibiting an explicit
graph $G$ with chromatic number $5$, such that $\Hom(G,K_5)$ is 
$0$-connected.
For the purpose of the present note, no background is needed beyond elementary
graph theory, and we refer the interested reader to~\cite{BaKo03,BaKo04}.

In the case under consideration here, $k=0$ and the graph $G$ has chromatic 
number $\chi$. (In the specific counterexample
we show, $\chi = 5$, but many other examples can be exhibited
for other values of $\chi$). 
The vertex set $W$ of the complex $Hom(G,K_{\chi})$
consists of all $\chi$-colorings of $G$. Whether this complex
is $0$-connected, is determined by its $1$-skeleton. This is a graph $G'$
on vertex set $W$, where
two $\chi$-colorings of $G$ are adjacent in $G'$ if they differ
on exactly one vertex. The complex is $0$-connected iff
the graph $G'$ is connected.
We present here a graph $G$ of chromatic
number $\chi=5$, for which the graph $G'$ is connected,
contradicting the above conjecture.

Let $G$ be the following graph with $9$ vertices and $22$ edges:
\begin{center}
\epsfig{file=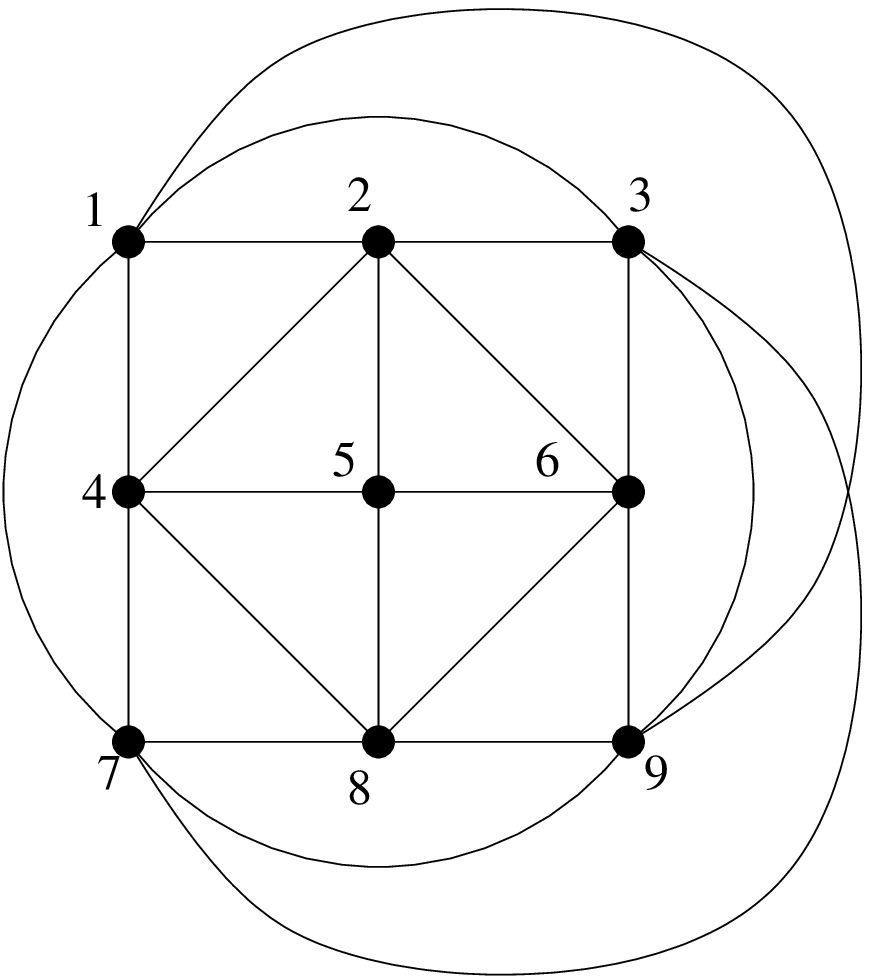,height=5cm}
\end{center}

It is easy to check that $\chi(G)=5$. We enumerate the 5-colorings of $G$
as follows. Since the four corners of the square form a clique, 
their colors must be distinct. 
Assume they are colored 1, 2, 3 and 4 in clockwise order starting at vertex
number 1 (this is one of the $5!$ ways to color the corners).
We would like to enumerate the ways to complete the following square:

{\gr 1-2 --- 4-3}

If the center is 5, then there are exactly two colorings:

$s_{1234}=${\gr 132 254 413}, $t_{1234}=${\gr 142 351 423}.

Otherwise, there are exactly 32 colorings. 
To see this, observe that if the center is different from 5, then 
the occurrences of 5 are either restricted to the central row or the central 
column but not to both.
This allows us to partition these colorings to $h$-type and $v$-type 
respectively.
Therefore, the colorings where the center is different from 5 are the possible
completions of the following eight squares (the plus sign marks vertices that
must be colored 5).

$h_{1234,a}=${\gr 132 +-- 423}, 
$h_{1234,b}=${\gr 142 --+ 413}, 
$v_{1234,a}=${\gr 1+2 3-4 4-3}, 
$v_{1234,b}=${\gr 1-2 2-1 4+3},
		       	    
$h_{1234,c}=${\gr 132 --- 413}, 
$h_{1234,d}=${\gr 142 --- 423}, 
$v_{1234,c}=${\gr 1-2 3-1 4-3}, 
$v_{1234,d}=${\gr 1-2 2-4 4-3}.

One can verify that each of these eight squares have four possible 
completions.
For example, 
the central row of $h_{1234,a}$ colorings must be in $\{514, 515, 545, 541\}$
and 
the central row of $h_{1234,c}$ colorings must be in $\{245, 545, 525, 524\}$.
Since all resulting colorings are distinct, we deduce that the total 
number of 5-colorings is $5! \cdot (32+2) = 4080$.

To prove that the graph $G'$ is connected, first note that the 4 colorings 
represented by each of the eight squares above, are connected. Indeed, 
the restriction of $G'$ to each of these four vertex sets is a path.
For example:

{\gr 132 514 423} $\rightarrow$ {\gr 132 515 423} $\rightarrow$ 
{\gr 132 545 423} $\rightarrow$ {\gr 132 541 423}.

\ignore{
As shown by the following paths:

{\gr 132 514 423} $\rightarrow$ {\gr 132 514 425} 
$\rightarrow$
{\gr 132 514 435} $\rightarrow$ {\gr 132 524 435} 
$\rightarrow$      	       
{\gr 132 524 415} $\rightarrow$ {\gr 132 524 413},
				       	       
{\gr 132 541 423} $\rightarrow$ {\gr 135 541 423}
$\rightarrow$
{\gr 125 541 423} $\rightarrow$ {\gr 125 531 423}
$\rightarrow$       	       
{\gr 145 531 423} $\rightarrow$ {\gr 142 531 423},

$h_{1234,a}$ is connected to $h_{1234,c}$ and $h_{1234,d}$. By symmetry,
we get that $h_{1234,*}$ are connected, as well as $v_{1234,*}$. 
}

We denote an edge in the graph $G'$ as a coloring, where one of the vertices
is given two colors. Using this notation, consider the following edges:

{\gr 132 514 42{3,5}}, {\gr 132 524 41{5,3}};
{\gr 13{2,5} 541 423}, {\gr 14{5,2} 531 423};
{\gr {1,5}42 235 413}, {\gr {5,1}32 245 413}.

These edges connect 
$(h_{1234,a}, v_{1254,a})$, $(v_{1254,a}, h_{1234,c})$;
$(h_{1234,a}, v_{1534,b})$, $(v_{1534,b}, h_{1234,d})$; 
and $(h_{1234,b},$ $v_{5234,b})$, $(v_{5234,b}, h_{1234,c})$ respectively.
Therefore, $h_{1234,a}, h_{1234,b}, h_{1234,c}, h_{1234,d}$ are connected, 
and we denote this connected set of colorings by $h_{1234}$. 
By symmetry, the same applies to the $v$-type colorings. 
Therefore, if the four corners are colored $a$, $b$, $c$ and $d$ (in clockwise
order starting at vertex 1) then we have four connected sets of colorings
$h_{abcd}, v_{abcd}, s_{abcd}$ and $t_{abcd}$.
This reduces the number of connected components to $5! \cdot 4 = 480$.

Next, consider the following edges:

{\gr 13{2,5} 514 423}, {\gr 132 514 42{3,5}}, 
{\gr {1,5}42 235 413}, {\gr 142 235 {4,5}13}.

These edges connect $h_{1234}$ to $v_{1534}$, $v_{1254}$, $v_{5234}$ and 
$v_{1235}$.
By symmetry, $v_{1234}$ is connected to 
$h_{1534}$, $h_{1254}$, $h_{5234}$ and $h_{1235}$.
Therefore, we can get from any $h_{abcd}$ either to $h_{1234}$ or $v_{1234}$.
Since $s_{1234}$ is connected to $v_{5234}$ and $h_{1534}$,
we have the path
$v_{5234}, s_{1234}, h_{1534}, v_{1234}, h_{5234}$, and hence that all the 
$h$ and $v$ colorings are connected.
Finally, we are done, since $t_{1234}$ is connected to $h_{5234}$. 
Therefore, the 5-colorings graph $G'$ is connected as claimed.

\bibliography{counter}
\end{document}